\def\P{\mathcal{P}}
\def\A{\mathcal{A}}
\def\D{\mathcal{D}}

\def\fil{\mathop{\mbox{Fil}}}
\def\cof{\mbox{cof}}
\def\cxa{C_\a^\x}

\def\la{\langle}
\def\ra{\rangle}
\def\z{\zeta}
\def\l{\lambda}

\def\min{{\rm min}}

\def\qed{$\Box$\medskip}

\def\z{\zeta}

\def\otp{{\rm otp}}
\def\cof{{\rm cof}}
\def\l{\lambda}
\newtheorem{theorem}{Theorem} \newtheorem{lemma}[theorem]{Lemma}
\newtheorem{remark}[theorem]{Remark}

\newtheorem{proposition}[theorem]{Proposition}
\newtheorem{observation}[theorem]{Observation}

\newtheorem{definition}[theorem]{Definition}

\def\a{\alpha}
\def\b{\beta}
\def\g{\gamma}

\def\e{\epsilon}
\def\x{\xi}
\def\d{\delta}
\def\p{\phi}
\def\o{\omega}
\def\k{\kappa}
\def\l{\lambda}
\def\sfin{\square^{\mbox{\tiny\it fin}}_{\l,D}}
\def\sfing{\square^{\g}_{\l,D}}
\def\sfinl{\square^{\l}_{\l,D}}
\def\sfinll{\square^{<\l}_{\l,D}}
\documentclass[12pt]{article}
\usepackage{amsfonts}
\usepackage{latexsym}
\usepackage{amssymb,amsmath,graphicx,psfrag}
\def\mm{\mathcal{M}}\def\mn{\mathcal{N}}
\def\II{II}\def\I{I}
\def\efg{EF}
\def\ef{Ehrenfeucht-Fra\"\i ss\'e}
\author{Juliette Kennedy\thanks{Research partially
supported by grant 251557 of the Academy of Finland.}
\\
Department of Mathematics and Statistics\\ University of Helsinki,
Finland\\
 \and Saharon Shelah \thanks{The second author would like
to thank the Israel Science Foundation for
      partial support of this research (Grant no. 1053/11). Publication 1011.}\\
Institute of Mathematics\\ Hebrew University, Jerusalem,
Israel\\
Rutgers University, New Jersey, USA\\
\and Jouko V\"a\"an\"anen\thanks{Research partially
supported by grant 251557 of the Academy of Finland.}\\
Department of Mathematics and Statistics\\ University of Helsinki,
Finland\\
ILLC, University of Amsterdam\\
Amsterdam, Netherlands }

\title{Regular Ultrapowers at Regular Cardinals}

\begin{document}

\maketitle

\begin{abstract}
In earlier work of the second and third author  the equivalence of a
finite square principle $\sfin$  with various model theoretic
properties of structures of size $\l$ and  regular ultrafilters was
established.
In this paper we investigate the principle $\sfin$, and thereby the above model theoretic properties, at a regular cardinal. By  Chang's Two-Cardinal Theorem, $\sfin$ holds at regular cardinals for all regular filters $D$ if we assume GCH. In this paper we prove in $ZFC$ that for certain regular filters that we call {\em doubly$^+$ regular}, $\sfin$ holds at regular cardinals, with no assumption about GCH. Thus we get new positive answers in $ZFC$ to  Open Problems 18 and 19  in the book {\em Model Theory} by Chang and Keisler.
\end{abstract}

\section{Introduction}\label{sec1}

In \cite{MR1926605} and \cite{MR2135667}  the equivalence of the following
finite square principle $\sfin$  with various model theoretic
properties of regular reduced powers of models was
established:

\begin{description}
\item[\(\sfin:\)]  $D$ is a filter on a cardinal $\lambda$ and
there exist finite sets \(\cxa\) and integers $n_\x$ for each
$\a<\lambda^+$ and $\x<\lambda$ such that for each \(\x,\a\)
\begin{description}

\item[(i)] \(\cxa\subseteq\a+1\)

\item[(ii)] If $B\subset\l^+$ is a finite set of ordinals and
$\a<\lambda^+$ is such that $B\subseteq\a+1$, then
$\{\x:B\subseteq \cxa\}\in D$

\item[(iii)] \(\b\in\cxa\) implies \(C_\b^\x=\cxa\cap(\b+1)\)

\item[(iv)] \(|\cxa|< n_\x\)

\end{description}
\end{description}

The model theoretic properties were the following: Firstly, if $D$
is an ultrafilter, then $\sfin$ is equivalent to $\mm^{\lambda}/D$
being $\lambda^{++}$-universal for each model $\mm$ in a vocabulary
of size $\le\l$. To formulate the second model theoretic property, let us say that two models are {\bf $EF_\a$-equivalent} if the second player (i.e. the ``isomorphism'' player) has a winning strategy in the \ef\ game of length $\a$ on the two models\footnote{The usual elementary equivalence in a finite relational vocabulary is thus $EF_n$-equivalence for all $n<\o$, and $L_{\infty\omega}$-equivalence is the same as $EF_\omega$-equivalence. For models of cardinality $\le \kappa$, $EF_\kappa$-equivalence is equivalent to isomorphism.}. Now    $\sfin$ is
equivalent to $\mm^\l/D$ and $\mn^\l/D$ being $EF_{\l^+}$-equivalent for any
elementarily equivalent models $\mm$ and $\mn$ (w.l.o.g. of cardinality $\le \l^+$) in a
vocabulary of size $\le\l$. The existence of such ultrafilters and
models is related to Open Problems 18 and 19  in the Chang-Keisler model theory book \cite{CK}.

The consistency of the failure of $\sfin$ for a regular filter at  a singular strong
limit cardinal $\l$ was proved  in \cite{MR2135667} relative to
the consistency of a supercompact cardinal. In \cite{MR2444269} this was improved to 
the failure of $\sfin$ for a regular {\em ultra}filter $D$ at a singular strong
limit cardinal $\l$  relative to
the consistency of a strongly compact cardinal. 
 The failure of $\sfin$ for an
ultrafilter implies the failure of $\lambda^{++}$-universality of
$\mm^{\lambda}/D$ for some $\mm$, as well as the failure of isomorphism of
some regular ultrapowers $\mm^\l/D$ and $\mn^\l/D$.  Thus \cite{MR2444269}
answered negatively the following problems listed in \cite{CK} modulo large
cardinal assumptions:
\begin{description}
\item[Problem 18 (\cite{CK})] Let $|M|,|N|, |{L}|\le\alpha$ and let $D$ be a regular ultrafilter over $\alpha$. If $\mm\equiv\mn$, then $\prod_D\mm\cong\prod_D\mn$.
\item[Problem 19 (\cite{CK})] If $D$ is a regular ultrafilter of $\alpha$, then for all infinite $\mm$, $\prod_D\mm$ is $\alpha^{++}$-universal.
\end{description} The use of large cardinals is justified by
 \cite{MR1926605}, \cite{MR2135667} and \cite{269} as the failure
of $\sfin$ for singular strong limit $\l$ implies the failure of
$\square_\lambda$, which implies the consistency of large
cardinals.

In this paper we investigate the principle $\sfin$, and thereby the above model theoretic problems, at a {\em regular} cardinal. The following result is proved in \cite{MR1147741}: Assume $\k$ is regular and $\l^{<\k}=\l$. Suppose $\mm$ and $\mn$ are structures for a finite vocabulary such that  $\mm$ and $\mn$ are $EF_\a$-equivalent for each $\a<\k$. Suppose $D$ is a filter on $\x\times\l$, $\x\le\l$, extending $F'\times F$, where $F'$ is a $\k$-descendingly incomplete filter on $\x$ and $F$ is a $\k$-semigood filter on $\l$ (the concept is defined in \cite{MR1147741}). Then 
$\mm^\l/D$ and $\mn^\l/D$ are $EF_{\l^+}$-equivalent. For $\k=\omega$ this, combined with the existence proof of semigood filters in \cite{MR1147741}, yields filters $D$
with $\sfin$.

The structure of the paper is the following: In Section~\ref{sec2} we prove weaker versions of $\sfin$ in the case where the filter $D$ extends the club filter on $\l$. Naturally this case is in spirit  quite far from the case of regular $D$, which is our prime interest. However, this result is useful in the sequel. Note that there are many regular (ultra)filters extending the club filter.  In Section~\ref{sec3}  we define the concept of {\em doubly$^+$ regular} filter and show that such filters $D$ on regular $\l>\aleph_0$ satisfy $\sfin$.  Thus we get new positive answers in $ZFC$ to the above Problem 18 (with isomorphism replaced, in the absence of $2^\l=\l^+$, by $EF_{\l^+}$-equivalence) and the above Problem 19. In Section~\ref{sec4} we prove results to the effect that not all regular filters are doubly regular. In Section ~\ref{sec5} we compare our concept of double regularity to Keisler's concept of goodness of a filter. In Section~\ref{sec6} we present some open questions.



\section{Filters extending the  club filter}\label{sec2}

We can get provable cases  of a weaker form  of $\sfin$, when $D$ extends the club filter.  This will prove useful in the next section, where we will use Theorem~\ref{a2} in the proof of Theorem~\ref{mainpoint}. The original $\sfin$ is equivalent to reduced powers of elementarily equivalent models of cardinality $\l$ being $EF_{\l^+}$-equivalent. The weaker form which we shall prove below will   give the $EF_{\l^+}$-equivalence of reduced powers of  models of power $\l$ that are not just elementarily equivalent but even $EF_\l$-equivalent.

\begin{theorem}\label{a2}
Suppose 
\begin{description}

\item[(a)] $\l$ is regular $>\aleph_0$, 
\item[(b)] $D$ is  a filter on $\l$.



\item[(c)] $D$ extends the club filter.

\end{description}
If $\mm$ and $\mn$ are $EF_\l$-equivalent, then $\mm^\l/D$ and $\mn^\l/D$ are $EF_{\l^+}$-equivalent.
\end{theorem}

\noindent{\bf Proof. }
If $\a<\l^+$, $\l$ regular, let $\{u^i_\a:i<\l\}$ be a continuously increasing sequence of subsets of $\a$ such that $|u^i_\a|<\l$ for all $i<\l$ and $\a=\bigcup_{i<\l}u^i_\a$. Let 
\begin{equation}\label{da}
D_\a=\{i<\l:\forall\b\in u^i_\a(u^i_\b=u^i_\a\cap\b)\}.\end{equation}
It is easy to see that $D_\a$ is a club of $\l$ (recall that $\l$ is regular). 

Now we can proceed,  as in \cite{MR1926605} to prove that if $M$ and $N$ are $EF_\l$-equivalent, then $M^\l/D$ and $N^\l/D$ are $EF_{\l^+}$-equivalent:

 Let $L$ be a finite vocabulary and for each
$i<\l$ let $\mm_i$ and $\mn_i$ be $EF_\l$-equivalent
$L$-structures. We show that \II\ has a winning strategy in the game
$\efg_{\l^+}$ on the models $\mm=\prod_D\mm_i$ and $\mn=\prod_D\mn_i$.

 The crucial
idea of the proof is the following: When  the \ef\ game
$\efg_{\l^+}(\mm,\mn)$ is played, the players are actually playing
$\l$ \ef\ games simultaneously, namely the games
$\efg_{\l}(\mm_i,\mn_i)$, $i<\l$. 

 For each $i<\l$ let $\sigma_i$ be a winning
strategy for \II\ in the game $\efg_{\l}$ on the models
$\mm_{i}$ and $\mn_{i}$.
A {\em good} position is  a sequence $\langle
(f_\b,g_\b):\b<\a\rangle$ for some $\a<\l^+$, together with a club $C\subseteq D_\a$, such that for all $\b<\a$ we
have $f_\b\in\prod_i M_{i}$, $g_\b\in\prod_i N_{i}$,
and if $i\in C$, then 
$$\langle (f_\eta(i),g_\eta(i)):\eta\in u^i_\a\rangle$$ is a play according to $\sigma_i$ on the models
$\mm_{i}$ and $\mn_{i}$. In a good position the
equivalence classes of the functions \(f_\b\) and \(g_\b\)
determine a partial isomorphism of the reduced products: Suppose \(\a\) rounds have been played
and we are in a good position. Let $\p_\g(
[f_{\b_1}],\ldots,[f_{\b_k}])$ be an atomic formula holding in
$\prod_i \mm_{i}/D$, where \(\b_1<\ldots<\b_k<\a\), and let
\(A=\{i\in D_\a:\{\b_1,\ldots,\b_k\}\subseteq u_{\alpha}^i\}\). By assumption, $A\in D$. Since also \(B=\{i<\l:{\mm_{i}}\models
\p_\g(f_{\b_1}(i),\ldots,f_{\b_k}(i))\}\in D\), we  have
\(A\cap B\in D\). For \(i\in A\cap
B\) we have $\b_1,\ldots,\b_{k}\in
u_{\a}^i$, hence $$u^i_{\b_j}=u^i_\a\cap\b_j.$$ Since we are in a good position,
\(\langle(f_\eta(i), g_\eta(i)):\eta\in u_{\a}^i\rangle\) is
a play according to winning strategy \(\sigma_i\). Hence \(\langle(f_\epsilon(\x),
g_\epsilon(\x)):\epsilon\in u_{\a}^i\rangle\) determines a
partial isomorphism of the structures \(\mm_{i}\) and \(\mn_{i}\). Since this was
the case for all $i\in A\cap B \in D$, we get $\prod_\e
\mn_{\e}/D\models \p_\g([g_{\b_1}],\ldots,[g_{\b_k}])$.
The strategy of \II\ is to keep the position of the game {\em
good} and thereby win the game. So suppose \(\b\) rounds have been
played and \II\ has been able to keep the position {\em
good}. 
Then for all $\g<\b$ there is a club $C_\g\subseteq D_\g$ such that for $i\in C_\g$, $\langle(f_\eta(i),
g_\eta(i)):\eta\in u^i_{\g}\rangle$ is a play according to
$\sigma_i$.
\medskip

\noindent {\bf Case 1:} $\b=\cup\b$. Let $C=\bigcap_{\g<\b}C_\g$. Since $\l$
is regular, this is still a club. We show that 
 $\langle (f_\g,g_\g):\g<\b\rangle$ is
good. Let \(i\in C\). Let us look at $\langle
(f_\eta(i),g_\eta(i)):\eta\in u_{\b}^i\rangle$. Since $i\in D_\b$, every initial segment of this play is a play according to $\sigma_i$.
Hence so is the entire play $\langle (f_\g,g_\g):\g<\b\rangle$. We have shown
that \II\ can maintain a good position. 

\noindent {\bf Case 2:} $\b=\d+1$. Let $C\subseteq\bigcap_{\g\le \b}C_\g$ such that $\d\in u^i_\b$ for $i\in C$. Now suppose \I\ plays \(f_\d\). We show that \II\ can
play \(g_\d\) so that $\langle (f_\g,g_\g):\g<\b\rangle$ remains
good. Let \(i\in C\). Let us look at $\langle
(f_\eta(i),g_\eta(i)):\eta\in u_{\d}^i\rangle$.
This is a play according to the strategy \(\sigma_i\). Since $i\in D_\b$ and $\d\in u^i_\b$,  
\(u_{\d}^i=u_\b^i\cap\d\), so after the moves $\langle
(f_\eta(i),g_\eta(i)):\eta\in u_{\d}^i\rangle$ \II\ can play one more move in \(EF_{\l}\)
on \(\mm_{i}\) and \(\mn_{i}\) with \I\ playing the element
\(f_\d(i)\). Let \(g_\d(i)\) be the answer of \II\ in this game
according to \(\sigma_i\). The values \(g_\d(i)\),
\(i\in C\), constitute the function \(g_\d\) mod $D$. We have shown
that \II\ can maintain a good position. 
\qed

\medskip

We do not know whether the conditions (a)-(c) of Theorem~\ref{a2} are necessary for the conclusion.

\begin{remark} We point out some variants of Theorem~\ref{a2}:
\begin{enumerate}

\item 
We can define a version $\sfing$ of $\sfin$ which is equivalent to: ``If $\mm$ and $\mn$ are $EF_\g$-equivalent, then $\mm^\l/D$ and $\mn^\l/D$ are $EF_{\l^+}$-equivalent": 
\begin{description}
\item[\(\sfing:\)]  $D$ is a filter on a cardinal $\lambda$ and
there exist finite sets \(\cxa\) and ordinals $\g_\x<\g$ for each
$\a<\lambda^+$ and $\x<\lambda$ such that for each \(\x,\a\)
\begin{description}

\item[(i)] \(\cxa\subseteq\a+1\)

\item[(ii)] If $B\subset\l^+$ is a set of ordinals with $\otp(B)<\g$ and
$\a<\lambda^+$ is such that $B\subseteq\a+1$, then
$\{\x:B\subseteq \cxa\}\in D$.

\item[(iii)] \(\b\in\cxa\) implies \(C_\b^\x=\cxa\cap(\b+1)\).

\item[(iv)] \(\otp(\cxa)< \g_\x\).

\end{description}
\end{description}

\noindent If clauses (a), (b) and (c) of Theorem~\ref{a2} are assumed, then $\sfinl$. 
\item 
 We can also define a version $\square^{<\d}_{\l,D}$ of $\sfin$ which is equivalent to ``If $\mm$ and $\mn$ are $EF_\g$-equivalent for all $\g<\d$, then $\mm^\l/D$ and $\mn^\l/D$ are $EF_{\l^+}$-equivalent". If clauses (a), (b) and (c)$^+$ of Theorem~\ref{a2} are assumed, then $\square^{<\l}_{\l,D}$ holds, where $(c)^+$ says that (c) holds and there are functions $f_\a$, $\a\le\l^+$, such that $\a<\b\le \l^+$ implies $\{i<\l:f_\a(i)<f_\b(i)\}\in D$ (For $D=$ the club filter this is the so called assumption of the existence of the $\l^+$'th canonical function, see e.g. \cite[p. 445]{jech}.) 
\item 
Note that $$\sfin\Rightarrow\sfing\Rightarrow\sfinll\Rightarrow\sfinl$$  
for $\g<\l$.

\item 
We get a variant of Theorem~\ref{a2} also by showing, assuming (a), (b) and (c), that $\prod_D\mm_i$ and $\prod_D\mn_i$ are $EF_{\l^+}$-equivalent, if for all $\b<\l$: $$\{i<\l:\mbox{$\mm_i$ and $\mn_i$ are $EF_\b$-equivalent}\}\in D.$$

\item We can weaken clause (c) of the theorem to the assumption that $D$ is  unreasonable (\cite{MR2283626}) in the following sense: There is a partition $\{w_i:i<\l\}$ of $\l$ such that $\bigcup_{i\in E}w_i\in D$ for every club $E$ of $\l$.

\end{enumerate}

\end{remark}

Û

\section{Doubly regular filters}\label{sec3}

We define the concept of a {\em doubly regular} filter, give examples of such on regular cardinals, and prove that $\sfin$ holds for such filters. Recall that a family of sets is a {\em regular family} if finite intersections of members of the family are non-empty, but all infinite intersections are empty, a filter  is called {\em $\mu$-regular} if it contains a regular family of size $\mu$, and a filter on $\lambda$ is called {\em regular} if it $\lambda$-regular. 

\begin{definition}\label{doubly} Suppose
 $D$ is a filter on a regular cardinal $\l$.
\begin{enumerate}

\item $D$
 is called {\em doubly regular}, if there are pairwise disjoint sets $u_i\subseteq\lambda$, $i<\l$, each of cardinality $\l$, and regular filters $D_i$ on $u_i$ such that for all $A\subseteq\l$: $$[\forall^{\infty}i<\l(A\cap u_i\in D_i)]\Rightarrow A\in D.$$(``\,$\forall^{\infty}i<\l$'' means ``for all but boundedly many $i$''.)

\item 
The filter $D$ is called {\em doubly${}^+$ regular} if the above holds with 
``$\, \forall^{\infty}i<\l$"
replaced by ``for a club of $i$". 
\end{enumerate}\end{definition}


Let us make some easy observations about doubly regular filters: 

\begin{observation}\label{obs}
\begin{enumerate}

\item  A doubly regular filter is necessarily regular: Let $\{A^\a_i:\a<\l\}$ be a regular family in $D_i$. Let $$B^\a=\bigcup_{i<\l}A^\a_i.$$ Then $\{B^\a:\a<\l\}$ is a regular family in $D$. We will show below (Theorem~\ref{last}) that the converse need not be true. 

\item A doubly$^+$ regular filter is always doubly regular.

\item 
It is easy to construct doubly$(^+)$ regular filters. Indeed, if the sets $u_i\subseteq\l, i<\l$, are disjoint, each of cardinality $\l$, $\l=\bigcup_i u_i$, and we have regular filters $D_i$ on $u_i$, then the set $\{A\subseteq \l: \forall^{\infty}i<\l(A\cap u_i\in D_i)\}$ is a doubly regular filter on $\l$, and 
the larger set $\{A\subseteq \l: \mbox{For a club of }i<\l(A\cap u_i\in D_i)\}$ is a doubly${}^+$ regular filter on $\l$. Both double regularity and double$^+$ regularity are closed under extensions of the filter, so we get also ultrafilter examples of both.

\end{enumerate}
\end{observation}

Here is the main point of doubly$^+$ regular filters, at least from the point of view  of this paper:

\begin{theorem}\label{mainpoint}
If $D$ is a doubly$^+$ regular filter on a regular cardinal $\l>\aleph_0$, then $\sfin$ holds.
\end{theorem}

\noindent{\bf Proof. }
Let the sets $u_i$ and the filters $D_i$ be as in Definition~\ref{doubly}. 
Let $D^*$ be the club filter  of $\l$, and $$D'=\{A\subseteq\l:\{i<\l:A\cap u_i\in D_i\}\in D^*\}.$$ We prove $\square^{\mbox{\tiny\it fin}}_{\l,D'}$. From this $\sfin$ follows, as $D'\subseteq D$.   It suffices to prove that if $\mm_\a$ and $\mn_\a$, $\a<\l$, are elementarily equivalent, with a vocabulary of size $\le\l$, then $\mm=\prod_{D'}\mm_\a$ and $\mn=\prod_{D'}\mn_\a$ are $EF_{\l^+}$-equivalent. Note that

\begin{description}
\item[(a)] $\mm\cong\prod_{i<\l}\mm^i/D^*$, where $\mm^i=\prod_{\a\in u_i}\mm_\a/D_i$.
\item[(b)] $\mn\cong\prod_{i<\l}\mn^i/D^*$, where $\mn^i=\prod_{\a\in u_i}\mn_\a/D_i$.
\end{description}
Since each $D_i$ is $\l$-regular, the models $\mm^i$ and $\mn^i$ are $EF_\l$-equivalent by \cite[Theorem VI.1.8]{MR1083551}. By Theorem~\ref{a2} the models $\mm$ and $\mn$ are now $EF_{\l^+}$-equivalent. \qed

\section{On regular but non-doubly regular filters}\label{sec4}

Non-regular uniform filters do not necessarily exist. If there is a non-regular uniform ultrafilter on $\omega_1$, then $V\ne L$ by \cite{prikry}, $0^\#$ exists by \cite{MR0419236}, and in fact $\omega_2$ is a limit of measurable cardinals in the Jensen-Dodd Core Model, by \cite{MR2000073}. We show that we can always construct a regular but non-doubly regular filter. In this sense double regularity is easier to avoid than regularity.

If $E$ is an equivalence relation on $\l$ we denote the set of all $E$-classes by $\l/E$, and the $E$-class of $i$ by $i/E$.

First an equivalent condition for double regularity, one that fits better our present purpose:

\begin{lemma}
A filter $D$ is doubly regular if and only if there is an equivalence relation $E$ of $\l$ and  $\bar{u}=\la u_\a:\a\in\l\ra$ such that:
\begin{description}

\item[(DR-a)] $\{u_\e:\e\sim_E i\}$ is a regular family  of subsets of $i/E$ for each $i<\l$.
\item[(DR-b)] If $S\subseteq\l$ and $|S|<\l$, then $\bigcup\{i/E:i\in S\}=\emptyset \mod D$,

\item[(DR-c)] $|i/E|=\l$ for all $i<\l$,

\item[(DR-d)] If $f$ is a function such that $dom(f)=\l/E$ and $f(i/E)\sim_E i$ for all $i\in \l/E$,  then $\bigcup_{i\in \l/E}u_{f(i)}\notin D$.

\end{description}
\end{lemma}

The proof is easy.



\begin{theorem}\label{last}
If $2^\l=\l^+$, then there is a regular ultrafilter on $\l$, which is not doubly regular.
\end{theorem}

\def\u{\bar{u}}

\noindent{\bf Proof. }
Let $\{B_\a: \a\in\l^+\}$ list $\P(\l)$. Let $\{(E^\a,\u_\a):\a<\lambda^+\}$ list potential candidates for double regularity i.e.  $E$ and 
${\u}=\la u_{\zeta}:\zeta<\l\ra$ such that $\{u_{\zeta}:\zeta<i/E\}$ is a regular family on $i/E$ for each $i<\l$. This is only place where we use $2^\lambda=\lambda^+$. 

We construct by induction sets $\D_\a$, $\a<\l^+$, such that the following conditions will hold:

\begin{description}

\item[(C-a)] $\D_\a\subseteq \P(\l)$ is $\subseteq$-continuously increasing.

\item[(C-b)] $|\D_\a|=\l$.

\item[(C-c)] $\D_\a$ is closed under finite intersections. We use $\fil(\D_\a)$ to denote the filter $\D_\a$ generates.

\item[(C-d)] $\D_0$ contains a regular family. (So necessarily, $u\in [\l]^{<\l}$ implies $u=\emptyset\mod D$.) 

\item[(C-e)] If $\alpha=2\beta+1$, then  $B_\beta\in\D_\a$ or $(\l\setminus B_\beta)\in\D_\a$.

\item[(C-f)] If $\alpha=2\beta+2$, then either there is $S\in[\l]^{<\l}$ such that $\bigcup_{\e\in S}\e/E_\b\ne\emptyset\mod\fil(\D_{\a})$, or, letting $\u_\b=\la u_{\b,\e}:\e<\l\ra$, 
there is $f$ such that $dom(f)=\l/E_\b$ and $f(i/E_\b)\sim_{E_\b} i$ for all $i\in \l/E_\b$,  then $\bigcup_{i\in \l/E_\b}u_{\b,f(i)}\in \D_a$.
\end{description}

Here is the construction:

\noindent {\bf Case 1:} $\a=0$. Let $E$ be a regular family on $\l$. (We can construct a regular family on $\l$ in the standard way: Let  $J$ be the set of finite subsets of $\l$. The family $\{\{X\in J:\beta\in X\}:\beta<\l\}$ is a regular family on $J$, and hence gives rise to one on $\l$.) We extend $E$ to  $\D_0$ by closing under finite intersections.

\medskip

\noindent {\bf Case 2:} $\a=2\beta+1$. We make a choice between $B_\beta\in\D_\a$ and $(\l\setminus B_\beta)\in\D_\a$ so that $\emptyset\notin\fil(\D_\a)$.
\medskip

\noindent {\bf Case 3:} $\a=2\beta+2$. Let $\{C^\a_l: l<\l\}$ list $\D_{2\b+1}$.
If there is $S\in[\l]^{<\l}$ such that $\bigcup_{\e\in S}\e/E_\b\ne\emptyset\mod\fil(\D_{2\b+1})$, we let $\D_{2\b+2}=\D_{2\b+1}$. So let us assume 
\begin{description}
\item[($\star$)] For all $S\in[\l]^{<\l}$ we have $\bigcup_{\e\in S}\e/E_\b=\emptyset\mod\fil(\D_{2\b+1})$.
\end{description}
\smallskip

We prove the following auxiliary:
\medskip

\noindent{\bf Subclaim:} There are $(\e_i,\g_i), i<\l$ such that 

\begin{center}
\begin{itemize}
\item[(a)] $\e_i\in\l\setminus\{\e_j:j<i\}$.

\item[(b)] $\g_i\sim_{E_\b} \e_i$. 

\item[(c)] $u_{\b,{\g_i}}\not\supseteq C^\a_i\cap \e_i/E_{\b}$.

\end{itemize} 

\end{center}
\medskip

 Let us first suppose the subclaim is true and we have such a sequence 
$(\e_i,\g_i), i<\l$. Choose $f$ by letting $f(\e_i)=\g_i$. So $\bigcup_{i\in \l/E_\b}u_{\b,f(i)}$ is a subset of $\l$, which includes no element of $\D_{2\b+1}$.
So we let $$\D_\a=\D_{2\b+1}\cup\{A\setminus\bigcup_{i\in \l/E_\b}u_{\b,f(i)}:A\in\D_{2\b+1}\}.$$ 
This is clearly closed under finite intersections and does not contain $\emptyset$ and every set in $\D_\a$ has cardinality $\l$.

Let us then prove the subclaim.  Let $i<\l$  and   $$W_1=\bigcup_{j<i}{\e_j}/E_\b.$$ By our assumption ($\star$), $W_1=\emptyset\mod\fil(\D_{2\b+1})$. Choose $\xi_i$ from the non-empty set $(\l\setminus W_1)\cap C_{\a,i}$. Then pick $\e_i$ so that $\xi_i\sim_{E_\b}{\e_i}$.
Finally, let 
$$W_2=\{\g<\l : \g\sim_{E_\b} {\e_i}\mbox{ and }\xi_i\in u_{\b,\g}\}.$$ 
Since ${\A}^\b$ is a regular family, the set $W_2$ is finite. So there is $\g_i\in u_{\b,{\e_i}}\setminus W_2$. This ends the construction of the sequence $(\e_i,\g_i), i<\l$, and thereby finishes the proof of the subclaim.
\medskip

\noindent{\bf Finishing the proof: } Now that we have constructed the sequence $\D_\a, \a<\l^+$, we can let $$D=\bigcup_{\a<\l^+}\D_\a.$$ This is an ultrafilter on $\l$. It is regular by (C-d). Now we can easily see that $D$ is not doubly regular: Suppose $E_\b$ and $\u_\b$ witnesses that $D$ is doubly regular. Let us look at the construction of $\D_{2\b+2}$. In the first case we assumed that there is $S\in[\l]^{<\l}$ with $\bigcup_{\e\in S}\e/E_\b\ne\emptyset\mod\fil(\D_{2\b+1})$. So $\bigcup_{\e\in S}\e/E_\b\ne\emptyset\mod D$, and (DR-b) is violated. In the second case we found $f$ such that $\bigcup_{i\in\l/E_\b}u_{\b,{f(i)}}=\emptyset\mod\fil(\D_\a)$. Hence $\bigcup_{i\in\l/E_\b}u_{\b,{f(i)}}=\emptyset\mod D$, and (DR-d) is violated. \qed

Note that  
    double$^+$ regularity of $D$ implies $\sfin$ on a regular cardinal $\l>\aleph_0$ (Theorem~\ref{mainpoint}), but in the light of the above Theorem, not conversely, as GCH implies $\sfin$ for regular $D$ and regular $\l$ (\cite[Lemma 4]{MR1926605}).


Theorem \ref{last} has the assumption $2^\l=\l^+$, which may fail for all $\l$. We shall present next a slightly different construction under a different assumption, one that is always satisfied by a multitude of cardinals $\l$.

\begin{theorem}\label{second}
Assume the following two conditions:\begin{description}

\item [(A1)]
 $\cof(\l)>\aleph_0$ or $\l>2^{\aleph_0}$. 
 \item [(A2)]
There is $\A\subseteq\P(\l)$ of cardinality $2^\l$ such that $|\{A\cap i:A\in\A\}|\le \l$ for all $i<\l$.
 
\end{description} Then there is a regular but not doubly regular filter on $\l$.
\end{theorem}

\noindent Note a family $\A$, as in (A2), always exists if $\l=2^{<\l}$. Hence condition (A2) can be replaced by $\l=\beth_\alpha$, $\alpha$ limit.
\medskip

\noindent{\bf Proof. }
Let $\la (E_\b,\u_\b):\b<2^\l\ra$ list all pairs where $E_\b$ is an equivalence relation on $\l$ and $\u^i_\b=\la u_{\b,\e}:\e\sim_{E_\b}i\ra$ is a regular family of subsets of $i/E_{\b}$ for each $i<\l$. Let $\{B_\a:\a<2^\l\}$ list $\P(\l)$.

We construct a sequence $(I_\a,\D_\a), \a<2^\l$ such that:  

\begin{enumerate}

\item 
$|I_\a|\le |\a|$, $I_\a\subseteq\P(\l)$, $(I_\a)$ is continuously increasing,

\item $\D_\a$ is the filter $\D[I_\a]=\{A\subseteq\l:\exists J\in
[I_\a]^{<\aleph_0}\exists S\in [\l]^{<\l}(\bigcap J\subseteq A\cup S)\},$

\item 
 $\D_{2\b+1}=\D_{2\b}\cup\{B_\b\}$ or $\D_{2\b+1}=\D_{2\b}\cup\{\l\setminus B_\b\}$

\item  $\D_{2\b+2}$ satisfies 

\begin{enumerate}
\item  There is some $W\in [\l]^{<\l}$ such that $\bigcup_{i\in W}i/E_\b\ne\emptyset$ mod $\D_{2\b+1}$, or

\item There is an $f$ such that $f(i/E_\b)\in i/E_\b$ 
for all $i$ and $\l\setminus\bigcup\{u_{\b,f(x)}:
x\in\l/E_\b\}\in I_\b$ or
\item $|\{X\in\l/E_\b : |X\cap B|=\l\}|<\l$ for some $B\in \D_{2\b+1}$.

\end{enumerate}

\end{enumerate}

The construction now follows:
Let us look at the case $\a=2\b+2$. If we cannot form $\D_\a$ as required, then:
\begin{description}

\item[(N1)] If $W\in [\l]^{<\l}$, then $\bigcup_{i\in W}i/E_\b=\emptyset$ mod $\D_{2\b+1}$.
\item[(N2)] If $f$ is a function such that $dom(f)=\l/E_\b$ and $f(i/E_\b)\sim_{E_\b}i$ for all $i<\l$, and $$A_{\b,f}=\bigcup\{u_{\b,f(x)}:x\in\l/E_\b\},$$then $\emptyset\in \D(I_{2\b+1}\cup \{\l\setminus A_{\b,f}\})$.
\item[(N3)] For $B\in \D_{2\b+1}$, $|\{X\in \l/E_\b:|X\cap B|=\l\}|=\l$.

\end{description}

We  derive a contradiction. This will ensure that $\D_\a$ can be found.  Let $\la x_{\b,i}:i<\l\ra$ list $\l/E_\b$. By our choice of $\A$, there are one-one functions $b_i:\{A\cap i:A\in \A \}\to x_{\b,i}$ for each $i<\l$. If $s\subseteq\l$, let $g_s$ be a function such that $dom(g_s)=\l/E_\b$ and $$g_s(x_{\b,i})=b_i(s\cap i)$$ so that $g_s(x_{\b,i})\in x_{\b,i}$. By (N2) there are 
$$J_{\b,s}\in [I_{2\b+1}]^{<\aleph_0}, W_{\b,s}\in[\l]^{<\l}$$ such that
$$\bigcap_{B\in J_{\b,s}}B\subseteq A_{\b,g_s}\cup W_{\b,s}.$$ Since $|\A|=2^\l$, there are $J_*\in [I_{2\b+1}]^{<\aleph_0}$ and $\mu<\l$ such that if 
$$\A_1=\{s\in\A:J_{\b,s}=J_*, |W_{\b,s}|=\mu \},$$ 
then $|\A_1|=2^\l$.  Let $B_*=\bigcap J_*\in \D_{2\b+1}$. By (N3), 
\begin{equation}\label{tama}
|\{j< \l:|x_{\b,j}\cap B_*|=\l\}|=\l.
\end{equation}
\bigskip

\noindent{\bf Claim:} There are $s_n\in\A_1$, $n<\omega$, and $i<\omega$ such that $s_n\cap i\ne s_m\cap i$ for all $n<m<\omega$.
\medskip

\noindent{\bf Case 1: $\cof(\l)>\aleph_0$.} Pick distinct $s_n\in\A_1$, $n<\o$. Since $\cof(\l)>\aleph_0$, there is $i<\l$ such that $s_n\cap i\ne s_m\cap i$ for all $n<m<\o$.  \bigskip

\noindent{\bf Case 2: $\cof(\l)=\aleph_0$, $\l>2^{\aleph_0}$.} Pick distinct $s_\x\in\A_1$, $\x<(2^{\aleph_0})^+$. Let $C\subseteq\l$ be cofinal, $|C|=\aleph_0$. Let $\chi:[(2^{\aleph_0})^+]^2\to C$ be defined by $\chi(\{\x,\z\})=\min\{c\in C:s_\x\cap c\ne s_\z\cap c\}$.  By the Erd\H{o}s-Rado Theorem $(2^{\aleph_0})^+\to (\aleph_1)^2_{\aleph_0}$ there is $i\in C$ and an uncountable $H\subseteq (2^{\aleph_0})^+$ such that $\chi\restriction [H]^2$ has constant value $i$.  
\medskip

The Claim is proved. By (\ref{tama}), there is $j>i$ such that $|B_*\cap x_{\b,j}|=\l$.
With the notation of (N2) $$A_{\b,g_{s_n}}\cap x_{\b,j}=u_{\b,b_j(s_n\cap j)}$$ and the sets $u_{\b,b_j(s_n\cap j)}$ are distinct because $b_j$ is one-one. By regularity, \begin{equation}\label{reg}\bigcap_n u_{\b,b_j(s_n\cap j)}=\emptyset.\end{equation}
Let $W=\bigcup\{W_{\b,s_n}:n<\o\}$. Clearly, $|W|=\mu$. Now $$B_*\cap x_{\b,j}\subseteq u_{\b,b_j(s_n\cap j)}\cup W.$$ This contradicts $|B_*\cap x_{\b,j}|=\l$, since  $|W|=\mu$   and (\ref{reg}) gives 
$$B_*\cap x_{\b,j}\subseteq \bigcap_n (u_{\b,b_j(s_n\cap j)}\cup W)=W.$$
\qed

If we start with a model of $GCH$, we can use  Easton forcing \cite{MR0269497} to obtain a model in which $2^{\l}$  is---for all regular $\l$---anything not ruled out by the conditions $\k\le \l\Rightarrow 2^{\k}\le 2^{\l}$ and $\cof(2^{\l})>\l$. In the arising forcing extension $V[G]$ the tree $({}^{<\l}2)^V$, $\l$ regular, has cardinality $\l$ and $2^\l$ branches. Hence we have in $V[G]$ a set $\A_\l$ of cardinality $2^\l$---for all regular $\l$---such that $\forall i<\l(|\{A\cap i : A\in \A_\l\}|\le\l)$, which is exactly the assumption (A2) of Theorem~\ref{second}.

\section{Good ultrafilters}\label{sec5}

Keisler
\cite{MR0166105} introduced the concept of $\k$-goodness of ultrafilters and proved that if $2^\l=\l^+$ and $D$ is a $\l^+$-good (i.e. good)  countably incomplete ultrafilter on $\l$, then $\prod_D\mm_i\cong\prod_D\mn_i$ for any models $\mm_i\equiv\mn_i$ of cardinality $\le \l^+$ in a vocabulary of cardinality $\le\l$.
This raises the question whether there is a connection between goodness and double regularity. It turns out that these concepts are independent of each other. 

\begin{proposition}Suppose $\l>\aleph_0$.
There is a doubly regular ultrafilter on $\l$ which is not good. If $2^\l=\l^+$, then there is a good countably incomplete ultrafilter on $\l$ which is not doubly regular.
\end{proposition}

\noindent{\bf Proof. } For the first claim, let $D_1$ be a doubly regular ultrafilter on $\l$ (exists by Observation~\ref{obs}) and $D_2$ a countably incomplete ultrafilter of $\omega$ which is not $\aleph_2$-good. (exists by \cite[5.1]{MR0166105}). Let $D=D_1\times D_2$. This is an ultrafilter on the set $\l\times\o$ of size $\l$. Since $D_2$ is not $\l^+$-good, neither is $D$ (\cite[VI.3.7]{MR1083551}). Double regularity is inherited from $D_1$ as follows: Suppose we have pairwise disjoint sets $u_i$, $i<\l$, on $\l$, each of cardinality $\l$, and regular filters $F_i$ on $u_i$ such that for all $A\subseteq \l$:$$[\forall^\infty i<\l(A\cap u_i\in F_i)]\rightarrow A\in D_1.$$ Let $G_i\subseteq F_i$ be a regular family on $u_i$. Let $u^*_i=u_i\times\o$ and $G^*_i=\{A\times\o:A\in G_i\}$. Let $F^*_i$ be the filter on $u^*_i$ generated by $\{A\times\o : A\in F_i\}$. Now $G^*_i$ is a regular family $\subseteq F^*_i$ and if $A\subseteq \l\times\o$, then  
$$[\forall^\infty i<\l(A\cap u^*_i\in F^*_i)]\rightarrow A\in D_1\times D_2.$$This ends the proof that $D$ is doubly regular. 

For the second claim we use a combination of the construction of the proof of Theorem~\ref{last} and Keisler's construction of a good ultrafilter in \cite[4.4]{MR0166105}. The construction of Keisler, as presented in \cite[Chapter 6, p. 387]{CK} proceeds in stages, generating a continuously increasing sequence $F_\a$, $\a<2^\l$, of filters such that the following condition holds (for unexplained terminology we refer to \cite[Chapter 6, p. 387]{CK}): For the first (in a fixed well-ordering) monotone $f:[\l]^{<\aleph_0}\to F_\a$ for which there is no additive extension $[\l]^{<\aleph_0}\to F_\a$, there is an additive extension $g:[\l]^{<\aleph_0}\to F_{\a+1}$. To make sure that such  $g$ and $F_{\a+1}$ always exist an auxiliary sequence is simultaneously defined, namely a descending sequence  $\Pi_\a$, $\a<2^\l$, of partitions of $\l$, starting from a carefully chose initial set $\Pi_0$ with $|\Pi_0|=2^\l$. There is no problem in interleaving the inductive construction of the filters $F_\a$ into the construction in the proof of Theorem~\ref{last}. The resulting ultrafilter is good but not doubly regular. \qed
 

\section{Concluding remarks}\label{sec6}

We proved that $\square^{\mbox{\tiny\it fin}}_{\lambda,D}$ holds if $\lambda$ is a regular cardinal and $D$ is a doubly regular  filter. This naturally raises the question whether $\square^{\mbox{\tiny\it fin}}_{\lambda,D}$ can fail at a regular cardinal for some regular, but not doubly regular, filter. We know it can fail at a singular cardinal \cite{MR2135667}. 
\bigskip

\noindent {\bf Conjecture 1:\ } Consistently,  $\square^{\mbox{\tiny\it fin}}_{\lambda,D}$ fails for some regular $\lambda>\omega$ and some regular filter $\lambda$ generated by $\lambda$ sets.
\medskip


\noindent {\bf Conjecture 2:\ } If $D$ is a regular ultrafilter on $\aleph_1$ such that $\neg \square^{\mbox{\tiny\it fin}}_{\aleph_1,D}$, then for any increasing continuous $\langle \alpha_i : i<\omega_1\rangle$ with $\alpha_i<\omega_1$, there is $A\in D$ such that $A\cap[\alpha_i,\alpha_{i+1})$ is finite for all $i<\omega_1$.
\medskip

\noindent Note that if  $$D=\{A\subseteq\omega_1:\forall^{\infty}i<\l(A\cap [\alpha_i,\alpha_{i+1})\in D_i\}),$$ $D_i$ ultrafilter on $[\alpha_i,\alpha_{i+1})$, then the answer to Conjecture 2 is positive. 
 This may indicate that looking for counterexamples for $\square^{\mbox{\tiny\it fin}}_{\aleph_1,D}$ can be hard.

\end{document}